# Using nodal coordinates as variables for the dimensional synthesis of mechanisms

V. García-Marina  · I. Fernández de Bustos · G. Urkullu · M. Abasolo

**Abstract** The method of the lower deformation energy has been successfully used for the synthesis of mechanisms for quite a while. It has shown to be a versatile, yet powerful method for assisting in the design of mechanisms. Until now, most of the implementations of this method used the dimensions of the mechanism as the synthesis variables, which has some advantages and some drawbacks. For example, the assembly configuration is not taken into account in the optimization process, and this means that the same initial configuration is used when computing the deformed positions in each synthesis point. This translates into a reduction of the total search space. A possible solution to this problem is the use of a set of initial coordinates as variables for the synthesis, which has been successfully applied to other methods. This also has some additional advantages, such as the fact that any generated mechanism can be assembled. Another advantage is that the fixed joint locations are also included in the optimization at no additional cost. But the change from dimensions to initial coordinates means a reformulation of the optimization problem when using derivatives if one wants them to be analytically derived. This paper tackles this reformulation, along with a proper comparison of the use of both alternatives using sequential quadratic programming methods. In order to do so, some examples are developed and studied.

V. García-Marina (✉)
Department of Mechanical Engineering, University School of Engineering of Vitoria-Gasteiz, University of the Basque Country, Nieves Cano 12, 01006 Vitoria-Gasteiz, Spain
e-mail: vanessa.garcia@ehu.eus

I. Fernández de Bustos · G. Urkullu · M. Abasolo
Department of Mechanical Engineering, Faculty of Engineering of Bilbao, University of the Basque Country, Alameda de Urquijo s/n, 48013 Bilbao, Spain
e-mail: igor.fernandezdebustos@ehu.eus

G. Urkullu
e-mail: gorka.urkullu@ehu.eus

M. Abasolo
e-mail: mikel.abasolo@ehu.eus

## 1 Introduction

The synthesis of mechanisms is a problem of high practical interest and, thus, it has been the scope of many research jobs. Synthesis problems can be classified in different types such as structural synthesis, geometrical synthesis, design synthesis, configuration synthesis, type synthesis, position synthesis, dimensional synthesis, kinetostatic synthesis, kinetic synthesis, kinematic synthesis, rectified synthesis, optimal synthesis, or probabilistic synthesis. The contributions presented in this paper will be centred in the kinematic dimensional synthesis, where the dimensions of the links of a mechanism are searched for while trying to fulfil certain position kinematic requisites defined in this case as synthesis points.

Actually, many methods have been used to accomplish the study of the synthesis of mechanisms and here a short resume will be presented. Some of these methods are heuristic and some fall into the numerical type of techniques. Between the first main group of techniques, there are the Genetic Algorithms [1–6], the Simulated Annealing [7–9], the Ant Colony Optimization [10, 11], the Particle Swarm Optimization [12–15], and some others like the Tabu Search [16, 17]. Among the second main group of techniques, there are the Sequential Quadratic Programming (SQP) methods where the most common ones are based on the method of Newton, or Quasi-Newton approaches. In order to introduce restrictions, if they are of linear nature, the Null Subspace method should be appropriate, and a good analysis of the different alternatives is exposed in [18]. If they are non-linear, the methods of the Penalty Function, of the Lagrange Multipliers, or the Augmented Lagrangian Function should be used. In the case of linear inequality restrictions, the methods of Karmarkar or the Primal-Dual should be adequate. Finally, for non-linear inequality restrictions, the methods of the Slack Variables or the Logarithmic Barrier Function can be used. [19–21].

A common way of classifying the types of synthesis problems is path generation, function generation, rigid solid guide and mixed synthesis. Path generation tries to obtain the best possible correlation between the path described by the joints of a mechanism during the solid rigid motion, together with some other previously specified path. Function generation studies the coordination or synchronization of the positions of the input and output links of a mechanism. Rigid solid guide is the part of the mechanism synthesis that studies the problem of locating a floating element (coupler) of a mechanism along a series of given positions. The mixed synthesis, in its turn, is a combination of some of the aforementioned types of synthesis in the same problem. In this paper a new alternative for a method for dimensional synthesis is presented, which has been under continuous development and accurate improvement for the last thirty years, since in 1982, in reference [22] for the first time the concept of the deformed position problem was presented. The main idea being to obtain the minimum energy position of the elements of a mechanism when one or more of its joints is obliged to fulfil certain geometrical restrictions out from the possible motions as rigid solid of the mechanism. The mechanism is considered composed of deformable elements with a linear elastic behaviour. Thus, the initial position problem was solved by means of the minimization of an energetic function, defined as the summation of the difference between deformed and undeformed squared lengths for each link in the mechanism. The same methodology was employed for the definition of the finite displacements problem, the deformed position problem and the static equilibrium problem. It was also suggested to solve the optimum synthesis based on these same ideas by summing the minimum deformation energy in each synthesis point. Later on, this idea was applied in [23], using the dimensions of the mechanism as variables. Exact derivatives were obtained for the deformed problem, but for the synthesis the length in each iteration was obtained via the arithmetical average of the deformed lengths of each of the deformed position problems. In 1989, the algorithm's convergence was improved by using approximate derivatives by means of finite differences in the synthesis instead of the arithmetical average [24]. The possibility of considering dimensional restrictions was also introduced by means of a stiffness that increases as long as the dimension gets nearer to the restricted value, which can be considered as a penalty function method. In 1993, in [25], the algorithm was improved by obtaining the exact first and second derivatives of every term. Furthermore, special elements with three joints were introduced to solve the low stiffness problems that appear whenever those joints are aligned. Here was also introduced a method to consider the fixed joints positions as variables. In order to do so, it is supposed that these are joined in the ends of a spring in direction of x and another spring in direction of y to a fixed point. In 2000, a preliminary study of the application of genetic algorithms to the synthesis of mechanisms and other mechanical problems is presented in [26]. Here it is demonstrated that the deformed position method is not very appropriate to be used together with the GA, due to the problem of the high aptitude of low stiffness mechanisms combined with the explorative nature of GAs. As a result, in 2004, another error function based on the minimum distance between the mechanism joint and the synthesis point [27] is applied, which happens to be valid to be applied with GAs.

In this paper the use of initial coordinates is explored for the synthesis of mechanisms using SQP and the deformation energy error function. The use of these kind of variables for the synthesis is not new, and has also been used using GAs and a distance error function with success, but the use of them in a deformation energy error function has not been yet studied. In this paper the relevant mathematical developments are presented and the analysis of several examples is exposed.

The motivation behind this change is that the use of dimensions does not include any information on the assembly configuration and, thus, the search space is somehow reduced. In the former formulation, an immutable set of coordinates were introduced by the user which were used to solve the position problem, and these coordinates had decisive influence on the deformed position problem solution.

This change has a main drawback though: an infinite set of initial coordinates define the same mechanism and, therefore, the optimization problem is always underdefined. This means that the optimization solver needs to be able to solve underdetermined systems.

The paper is organized as follows. First, a review of the deformed energy method is exposed. Afterwards, the choice of using initial coordinates as synthesis variables along with the deformed energy method is reasoned. Then, the energetic error function using initial coordinates is developed and the analytic expressions are presented, and the method for introducing boundary conditions is exposed. After that, some remarks on the optimization method are commented. Finally, some results are presented and some conclusions are driven.

## 2 The optimization of mechanisms using the deformed energy method

In order to better introduce the developments presented in this paper, a brief approach to the deformed energy method will be exposed here. To provide simplicity, this explanation will be reduced to mechanisms represented by truss elements joined by R-joints. In previous work, as shown in the introduction, the mechanism dimensions were defined by a vector whose components are those dimensions. The error function is tailored as follows: The mechanism under study is placed in an initial position, and expressed through the dimensions of the links, calculated by means of the nodal coordinates of all the joints of the mechanism. That is, the data for the problem are only the dimensions of the links of the mechanism, so by setting an initial position, an initial assembly configuration will be established, and this is usually not changed during the optimization process. Then the deformed position problem will be solved for each of the precision positions, by defining the nodal coordinates that give the optimum dimensions of all the links, usually slightly different for each of the precision points, so that the stored deformation energy in the whole mechanism is the minimum. The error function used for evaluating the fitness is the deformation energy, which is measured in each precision point and summed for all of the positions as it can be appreciated in Eq. (1).

$$F(L) = \sum_{i=1}^{P} \sum_{j=1}^{B} (L_j - l_{ji}(x_i))^2 \qquad (1)$$

This function represents the deformation energy of the B trusses of the mechanism supposing these are deformable to be able to reach the P desired synthesis positions. The $L_j$ are the non-deformed lengths and the $l_{ji}(x_i)$ are the lengths of the same trusses but now deformed in each of the i precision positions, and expressed through the nodal coordinates of the joints.

To give an example of what does the deformed problem consist on, a fourbar mechanism is going to be used. Let us place the mechanism in an initial position and let us define three precision positions for the coupler point as seen in Fig. 1. Here nodes A and B are fixed and node E (the coupler point) is the one to follow the prescribed path defined by precision points 0–1–2.

The mechanism solved for each of the three precision positions is shown in the next Fig. 2a–c, where it can be observed that the trusses are deformed with respect to the initial lengths and are different for each position.

In the synthesis error function, the deformed position problem is solved for each one of the precision positions, so that a set of coordinates is obtained for each of them, which, in term, define the deformed lengths $l_{ji}(x_i)$.

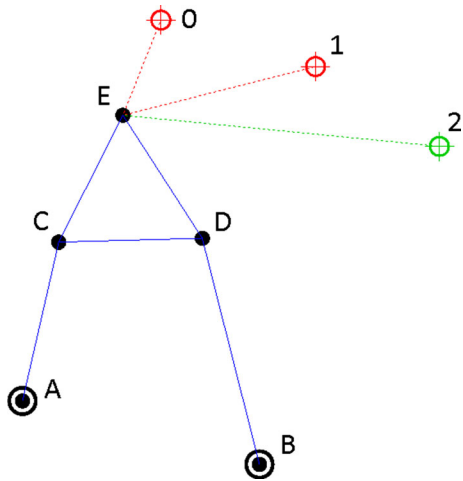

**Fig. 1** Fourbar in initial position and three precision positions

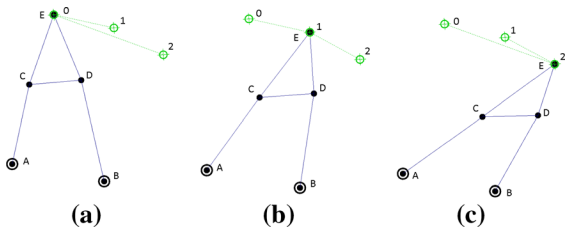

**Fig. 2** Fourbar in deformed positions **a** 0, **b** 1, and **c** 2

The optimization of the synthesis function has been approached in two ways. The uncoupled approach, which is based in discarding the effect of the modification of a dimension in the deformed position problem [23] and the coupled approach, which takes into account this effect [24, 25].

## 3 Reasons for the use of initial coordinates as parameters for the optimization of mechanisms

To keep the formulation simple, this paper will only include the definition of mechanisms composed by R-Joints and modeled by simple truss elements. No higher order elements or other joints will be considered, although the ideas exposed here can easily be generalized for developments such as those published in [25] or [28].

In order to clarify this point, it is first necessary to introduce how are usually defined the dimensions of a mechanism in the optimization process. Let us consider, for example a simple fourbar as that in Fig. 3. Again, fixed nodes are A and B. In this example instead of identifying nodes, we identify the links because in this case dimensional optimization parameters are the lengths.

The synthesis variable vector would in this case be defined as:

$$L^T = \{ L_0 \quad L_1 \quad L_2 \quad L_3 \quad L_4 \}$$

The use of the dimensions of the mechanism as parameters for the optimization is the most straightforward approach when performing dimensional optimal synthesis of mechanisms. It also has some additional advantages. For example, if one wants to introduce a determinate value for one dimension, this translates in this case as a linear restriction. This allows the use of simple methods such as the nullspace method (see, for example, [18] for a good review on cost-effective methods on introducing linear restrictions), without having to resort to Lagrange multipliers or other methods for the introduction of non linear restrictions. This also applies for interior point restrictions, where one can use Karmarkar or similar methods instead of Logarithm Barrier or others. But this does not come without drawbacks. One of the most important drawbacks is related to the assembly configuration. The dimensions of the mechanism do not define the assembly of the mechanism and, thus, one needs to somehow introduce the assembly configuration. For example, as can be seen in Fig. 4, mechanisms A and B have the same dimensions but is impossible to switch from one configuration to the other without dismantling the mechanism. Until now, this assembly configuration was defined by the user by introducing a set of initial coordinates for solving the deformed position problem. It is important to state that the assembly configuration introduced in this way will not always be maintained in the final result, but it has a determinant role in it. However, initial coordinates do define the assembly configuration of the mechanism, so making use of these as optimization parameters, the assembly configuration is being included in this optimization. Another drawback is derived from the fact that position of the fixed joints (those united to the fixed link) are not directly included in the optimization process, and one needs to modify the algorithm to optimize them. If one also wants to limit the space where those joints are to be optimized (restricted optimization), the algorithm gets quite complex.

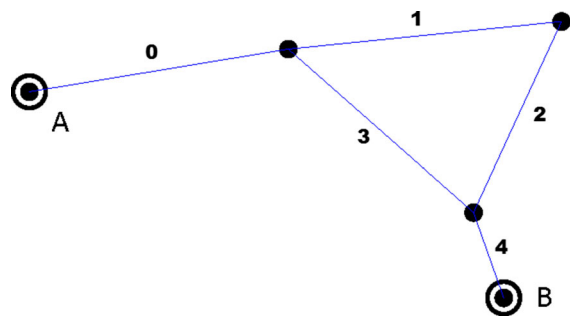

**Fig. 3** Topology of a fourbar, where optimization parameters are the dimensions

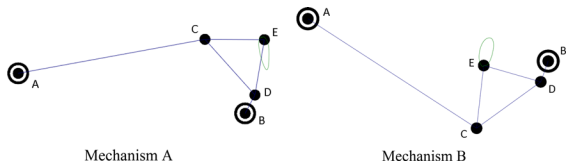

**Fig. 4** Two possible configurations with the same dimensions

The use of initial coordinates has been used in other synthesis methods as, for example, in [29] and it was also employed in [27] to tackle the problem of the assembly configuration when using genetic algorithms to optimize mechanisms, with quite a good result, and using an error function based in distances, instead of deformation energy. This was necessary due to the exploratory nature of genetic algorithms. Let us consider the examples exposed before. In the new formulation, the design vector is now composed of a set of initial coordinates as it shown in Fig. 5, where notations of each node are given.

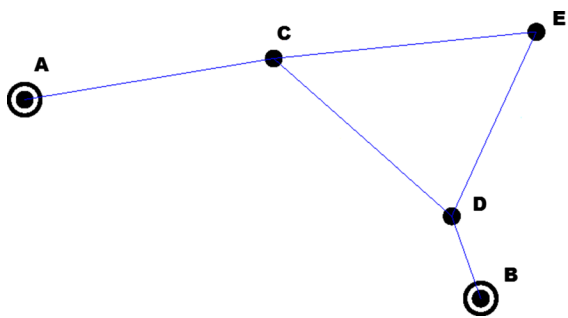

**Fig. 5** Topology of a fourbar, where optimization parameters are the nodal coordinates

In this case the synthesis vector of variables would have the form as follows:

$$x_0^T = \{\, x_{A0} \quad y_{A0} \quad x_{B0} \quad y_{B0} \quad x_{C0} \quad y_{C0} \quad x_{D0} \quad y_{D0} \quad x_{E0} \quad y_{E0} \,\}$$

In this paper this formulation is used, but considering an energy based error function and using SQP methods. These methods, although exploitative, can benefit from the other advantages of the formulation, namely the optimization of the fixed joints location and it is also adequate if one wants to perform mixed optimization combining exploratory and exploitative methods. It also comes with drawbacks, being the most relevant of them the redundancy of the solution vector, in the sense that different solution vectors may well represent the same mechanism. This leads to the need of using solvers capable of tackling with indefinite Hessian matrices.

## 4 Error function

Once the design vector is defined, one needs to specify the error function to be optimized. If one is to keep using a deformed energy error function, this error function must be rewritten to be expressed in terms of the initial coordinates. This is: instead of the Eq. (1) (here it is recalled that the presented formulation is limited to truss elements and R joints in spite of simplicity), one must use Eq. (2).

$$F(x_0) = \sum_{i=1}^{P} \sum_{j=1}^{B} \bigl(L_j(x_0) - l_{ji}(x_i)\bigr)^2 \quad (2)$$

where $P$ is the number of precision points; $B$ is the number of trusses defining the mechanism; $L_j$ is the dimension of the $j$-th truss (and, thus, the design vector, $L$, is of dimension $P$); $x_i$ are the set of coordinates which minimize the deformation energy of the mechanism for the requirements in the precision point $i$, and $l_{ji}$ is the length of the truss $j$ as defined by the set of coordinates $x_i$; this is, the deformed length of the $j$-th truss. It is important then to state that each of the $x_i$ vectors are obtained by an optimization process whose objective is to yield the lower deformation energy of the mechanism (considered as deformable) in the $i$-th precision point.

In Eq. (2), instead of being taken as optimization parameters, $L_j$ are defined by the optimization parameters $x_0$, which represent a set of initial coordinates. One of the strong points of Eq. (2) is that it can be expressed as the assembly of a finite element matrix composed of truss elements by using the form in Eq. (3):

$$F(x_0) = \sum_{i=1}^{P}\sum_{j=1}^{B}\left(L_j\{x_0\} - l_{ji}(x_i)\right)^2 = \sum_{j=1}^{B} F_{j0}^e(x_0, x_i) \quad (3)$$

This is quite convenient, because one can now perform operations on a per element basis.

## 5 Computing the derivatives

As exposed before, the computation of derivatives when using dimensions as variables has been performed in two different ways. The exact one, usually called coupled approach and an approximated one, called uncoupled. Both have their advantages and disadvantages. The coupled approach has second order convergence, but the derivation of the Hessian matrix is quite costly, while the uncoupled approach has lower rate of convergence but at a smaller iteration cost. The difference appears in the dependence of vectors $x_i$ on vector $L$. Taking the derivative of the expression in Eq. (1) with respect to $L_j$, one can write Eq. (4):

$$\frac{dF(L)}{dL_j} = \frac{\partial F(L)}{\partial L_j} + \sum_{i=1}^{P}\left(\frac{\partial F(L)}{\partial x_i}\right)^T \frac{\partial x_i}{\partial L_j} \quad (4)$$

The use of this full expression would be the so called coupled approach, while in the uncoupled approach one uses the expression in Eq. (5):

$$\frac{dF(L)}{dL_j} \approx \frac{\partial F(L)}{\partial L_j} \quad (5)$$

Here it will be demonstrated that, for this first derivative, one introduces no error when using Eq. (5) instead of (4). $x_i$ is the set of coordinates that delivers minimal deformation energy in the synthesis point $i$. Thus, one can write Eq. (6):

$$\frac{\partial \sum_{j=1}^{B}(L_j - l_{ji}(x_i))^2}{\partial x_i} = 0 \quad (6)$$

In the other hand:

$$\sum_{i=1}^{P}\left(\frac{\partial F(L)}{\partial x_i}\right)^T \frac{\partial x_i}{\partial L_j} = \sum_{i=1}^{P}\left(\frac{\partial \sum_{j=1}^{B}(L_j - l_{ji}(x_i))^2}{\partial x_i}\right)^T \frac{\partial x_i}{\partial L_j} = 0 \quad (7)$$

One cannot state the same for the second order derivatives. In this paper the uncoupled formulation will be applied, but using $x_0$ instead of $L$ a similar development to that presented in Eqs. (6) and (7) could be performed in this case, leading to similar conclusions, this is, the first derivatives are not affected by the use of coupled and uncoupled hypothesis. One can write the expression in (8).

$$g_j^e = \frac{\partial F_{j0}^e(x_{j0})}{\partial x_{j0}} = \begin{Bmatrix} \frac{\partial F(x_0)}{\partial x_{jk0}} \\ \frac{\partial F(x_0)}{\partial y_{jk0}} \\ \frac{\partial F(x_0)}{\partial x_{jl0}} \\ \frac{\partial F(x_0)}{\partial y_{jl0}} \end{Bmatrix}$$

$$= 2\sum_{i=1}^{P}\left(1 - \frac{l_{ji}(x_{ji})}{L_j(x_{j0})}\right)\begin{Bmatrix} x_{jk0} - x_{jl0} \\ y_{jk0} - y_{jl0} \\ x_{jl0} - x_{jk0} \\ y_{jl0} - y_{jk0} \end{Bmatrix}$$

$$= 2\left(P - \sum_{i=1}^{P}\frac{l_{ji}(x_{ji})}{L_j(x_{j0})}\right)\delta_{j0} \quad (8)$$

where:

$$x_{j0} = \begin{Bmatrix} x_{jk0} \\ y_{jk0} \\ x_{jl0} \\ y_{jl0} \end{Bmatrix}; \quad \delta_{j0} = \begin{Bmatrix} x_{jk0} - x_{jl0} \\ y_{jk0} - y_{jl0} \\ x_{jl0} - x_{jk0} \\ y_{jl0} - y_{jk0} \end{Bmatrix} \quad (9)$$

To reach to this equation and those which will follow, the identity in Eq. (10) is of importance:

$$\frac{\partial L_j(x_{j0})}{\partial x_{j0}} = \frac{1}{L_j(x_{j0})}\delta_{j0} \quad (10)$$

The second derivative (Hessian matrix) will have the form of expression in (11).

$$H_{j0}^e = \frac{\partial^2 F_{j0}^e(x_{j0})}{(\partial x_{j0})^2}$$

$$= \begin{bmatrix} \frac{\partial^2 F(x_{j0})}{(\partial x_{jk0})^2} & \frac{\partial^2 F(x_{j0})}{\partial x_{jk0}\partial y_{jk0}} & \frac{\partial^2 F(x_{j0})}{\partial x_{jk0}\partial x_{jl0}} & \frac{\partial^2 F(x_{j0})}{\partial x_{jk0}\partial y_{jl0}} \\ \frac{\partial^2 F(x_{j0})}{\partial y_{jk0}\partial x_{jk0}} & \frac{\partial^2 F(x_{j0})}{(\partial y_{jk0})^2} & \frac{\partial^2 F(x_{j0})}{\partial y_{jk0}\partial x_{jl0}} & \frac{\partial^2 F(x_{j0})}{\partial y_{jk0}\partial y_{jl0}} \\ \frac{\partial^2 F(x_{j0})}{\partial x_{jl0}\partial x_{jk0}} & \frac{\partial^2 F(x_{j0})}{\partial x_{jl0}\partial y_{jk0}} & \frac{\partial^2 F(x_{j0})}{(\partial x_{jl0})^2} & \frac{\partial^2 F(x_{j0})}{\partial x_{jl0}\partial y_{jl0}} \\ \frac{\partial^2 F(x_{j0})}{\partial y_{jl0}\partial x_{jk0}} & \frac{\partial^2 F(x_{j0})}{\partial y_{jl0}\partial y_{jk0}} & \frac{\partial^2 F(x_{j0})}{\partial y_{jl0}\partial x_{jl0}} & \frac{\partial^2 F(x_{j0})}{(\partial y_{jl0})^2} \end{bmatrix}$$

(11)

Taking into account the fact that:

$$\frac{\partial^2 L_j(x_{j0})}{(\partial x_{j0})^2} = \frac{1}{L_j(x_{j0})} \frac{\partial \delta_{j0}}{\partial x_{j0}} - \frac{1}{L_j^3(x_{j0})} \delta_{j0} \delta_{j0}^T \quad (12)$$

where:

$$\frac{\partial \delta_{j0}}{\partial x_{j0}} = \begin{bmatrix} 1 & 0 & -1 & 0 \\ 0 & 1 & 0 & -1 \\ -1 & 0 & 1 & 0 \\ 0 & -1 & 0 & 1 \end{bmatrix} \quad (13)$$

One can reach the expression in Eq. (14):

$$\begin{aligned} H_{j0}^e &= 2 \sum_{i=1}^{P} \left(1 - \frac{l_{ji}(x_{ji})}{L_j(x_{j0})}\right) \frac{\partial \delta_{j0}}{\partial x_{j0}} \\ &+ 2 \sum_{i=1}^{P} \delta_{j0} l_{ji}(x_{ji}) \frac{1}{L_j^2(x_{j0})} \left(\frac{\partial L_j(x_{j0})}{\partial x_{j0}}\right)^T \\ &= 2 \left(P - \sum_{i=1}^{P} \frac{l_{ji}(x_{ji})}{L_j(x_{j0})}\right) \frac{\partial \delta_{j0}}{\partial x_{j0}} \\ &+ 2 \sum_{i=1}^{P} \frac{l_{ji}(x_{ji})}{L_j^3(x_{j0})} \delta_{j0} \delta_{j0}^T \end{aligned} \quad (14)$$

This expressions allows one to obtain the full Hessian matrix and full gradient vector by means of expansion and assembly of the matrices obtained for each truss element.

## 6 Boundary conditions

The use of a set of nodal coordinates as parameters for the optimization also leads to a phenomenon to take into consideration: if any of these coordinates belong to a fixed joint in the mechanism, a change in these coordinates would also lead to a change in the deformed lengths obtained in each of the precision points. This is, $l_{ij}(x_{ji})$ would be affected and, thus, this phenomenon is to be considered in $g$ and $H$. If one considers Eq. (15):

$$\begin{aligned} \frac{dF(x_0)}{dx_0} &= \frac{\partial}{\partial x_0} \sum_{i=1}^{P} \sum_{j=1}^{B} \left(L_j(x_0) - l_{ji}(x_i)\right)^2 \\ &+ \sum_{i=1}^{P} \sum_{j=1}^{B} \frac{\partial x_i}{\partial x_0} \frac{\partial}{\partial x_i} \left(L_j(x_0) - l_{ji}(x_i)\right)^2 \end{aligned} \quad (15)$$

It is easy to find out that the effect of the boundary conditions in $g$ can be expressed as the summation of those components obtained in the previous section and an additional term. In order to use the same finite element approach described before, one can write Eq. (16):

$$\begin{aligned} \frac{dF(x_0)}{dx_0} &= \sum_{j=1}^{B} \frac{\partial}{\partial x_{j0}} \sum_{i=1}^{P} \left(L_j(x_{j0}) - l_{ji}(x_{ji})\right)^2 \\ &+ \sum_{j=1}^{B} \sum_{i=1}^{P} \left(\frac{\partial x_{ji}}{\partial x_{j0}}\right) \frac{\partial}{\partial x_{ji}} \left(L_j(x_{j0}) - l_{ji}(x_{ji})\right)^2 \\ &= \sum_{j=1}^{B} g_{j0}^e + \sum_{j=1}^{B} \sum_{i=1}^{P} g_{ji}^e \end{aligned} \quad (16)$$

So it can be considered the expression in Eq. (17), for a truss $j$, joining joints $k$ and $l$:

$$g_j^e = g_{j0}^e + \sum_{i=1}^{P} g_{ji}^e \quad (17)$$

Where $g_{j0}^e$ has already been obtained. Then one can write the equation (18):

$$\begin{aligned} g_{ji}^e &= \left(\frac{\partial x_{ji}}{\partial x_{j0}}\right) \frac{\partial}{\partial x_{ji}} \left(L_j(x_{j0}) - l_{ji}(x_{ji})\right)^2 \\ &= 2 \left(1 - \frac{L_j(x_{j0})}{l_{ji}(x_{ji})}\right) \frac{\partial x_{ji}}{\partial x_{j0}} \delta_{ji} \end{aligned} \quad (18)$$

Where expression in Eq. (19) are fulfilled.

$$x_{ji} = \begin{Bmatrix} x_{ki} \\ y_{ki} \\ x_{li} \\ y_{li} \end{Bmatrix}; \quad \delta_{ji} = \begin{Bmatrix} (x_{ki} - x_{li}) \\ (y_{ki} - y_{li}) \\ -(x_{ki} - x_{li}) \\ -(y_{ki} - y_{li}) \end{Bmatrix};$$

$$\frac{\partial x_{ji}}{\partial x_{j0}} = \begin{bmatrix} f_k & 0 & 0 & 0 \\ 0 & f_k & 0 & 0 \\ 0 & 0 & f_l & 0 \\ 0 & 0 & 0 & f_l \end{bmatrix} \quad (19)$$

$f_k$ equals 1 if node $k$ in truss $j$ is fixed and 0 if it is not; $f_l$ equals 1 if node $l$ in truss $j$ is fixed and 0 in the other case. It has also been used here the Eq. (20):

$$\frac{\partial l_{ji}(x_{ji})}{\partial (x_{ji})} = \frac{1}{l_{ji}(x_{ji})} \delta_{ji} \quad (20)$$

For $H_j^e$, one can write Eq. (21):

$$\begin{aligned} H_j^e &= \frac{dg_j^e}{dx_{j0}} = \frac{\partial g_{j0}^e}{\partial x_{j0}} + \sum_{i=1}^{P} \frac{\partial x_{ji}}{\partial x_{j0}} \frac{\partial g_{j0}^e}{\partial x_{ji}} \\ &+ \sum_{i=1}^{P} \frac{\partial g_{ji}^e}{\partial x_{j0}} + \sum_{i=1}^{P} \frac{\partial x_{ji}}{\partial x_{j0}} \frac{\partial g_{ji}^e}{\partial x_{ji}} \\ &= H_{j0}^e + \sum_{i=1}^{P} H_{j0i}^e + \sum_{i=1}^{P} H_{ji0}^e + \sum_{i=1}^{P} H_{jii}^e \end{aligned} \quad (21)$$

To derive the relevant matrices, the identity expressed in Eq. (22) is of extensive use:

$$\frac{\partial}{\partial x_{ji}}\left(1 - \frac{l_{ji}(x_{ji})}{L(x_{j0})}\right) = \frac{-1}{l_{ji}(x_{ji})} \delta_{ji} \quad (22)$$

$H_{j0}^e$ was attained in the previous Sect. 5. $H_{j0}^e$ is obtained from Eq. (23).

$$H_{j0ij}^e = \frac{\partial x_{ji}}{\partial x_{j0}} \frac{\partial g_{j0}^e}{\partial x_{ji}} = \frac{-2}{L(x_{j0})l(x_{ji})} \begin{bmatrix} f_k & 0 & 0 & 0 \\ 0 & f_k & 0 & 0 \\ 0 & 0 & f_l & 0 \\ 0 & 0 & 0 & f_l \end{bmatrix} \delta_{j0} \delta_{ji}^T \quad (23)$$

Which is gathered from expressions in Eq. (24):

$$\begin{aligned} \frac{\partial g_{j0}^e}{\partial x_{ji}} &= 2\frac{\partial}{\partial x_{ji}}\left(\left(1 - \frac{l_{ji}(x_{ji})}{L(x_{j0})}\right)\delta_{j0}\right) \\ &= 2\delta_{j0}\left(\frac{\partial}{\partial x_{ji}}\left(1 - \frac{l_{ji}(x_{ji})}{L(x_{j0})}\right)\right)^T \\ &\quad + 2\left(1 - \frac{l_{ji}(x_{ji})}{L(x_{j0})}\right)\frac{\partial}{\partial x_{ji}}\delta_{j0} \\ &= -2\delta_{j0}\frac{1}{L(x_{j0})l(x_{ji})}\delta_{ji}^T \end{aligned} \quad (24)$$

For $H_{jji0}^e$ it can be expressed as in Eq. (25) the following:

$$\begin{aligned} H_{jji0}^e &= \frac{\partial g_{ji}^e}{\partial x_{j0}} = 2\frac{\partial}{\partial x_{j0}}\left(\left(1 - \frac{L_j(x_{j0})}{l_{ji}(x_{ji})}\right)\frac{\partial x_{ji}}{\partial x_{j0}}\delta_{ji}\right) \\ &= 2\frac{\partial x_{ji}}{\partial x_{j0}}\delta_{ji}\left(\frac{\partial}{\partial x_{j0}}\left(1 - \frac{L_j(x_{j0})}{l_{ji}(x_{ji})}\right)\right)^T \\ &\quad + 2\left(1 - \frac{L_j(x_{j0})}{l_{ji}(x_{ji})}\right)\frac{\partial}{\partial x_{ji}}\left(\frac{\partial x_{ji}}{\partial x_{j0}}\delta_{ji}\right) \\ &= \frac{-2}{l_{ji}(x_{ji})L_j(x_{j0})}\frac{\partial x_{ji}}{\partial x_{j0}}\delta_{ji}\delta_{j0}^T \\ &\quad + 2\left(1 - \frac{L_j(x_{j0})}{l_{ji}(x_{ji})}\right)\frac{\partial}{\partial x_{ji}}\left(\frac{\partial x_{ji}}{\partial x_{j0}}\delta_{ji}\right) \end{aligned} \quad (25)$$

Where the last term in Eq. (25) can be expressed as shown in Eq. (26):

$$\frac{\partial}{\partial x_{j0}}\left(\frac{\partial x_{ji}}{\partial x_{j0}}\delta_{ji}\right) = \begin{bmatrix} f_k & 0 & -f_k f_l & 0 \\ 0 & f_k & 0 & -f_k f_l \\ -f_k f_l & 0 & f_l & 0 \\ 0 & -f_k f_l & 0 & f_l \end{bmatrix} \quad (26)$$

The last term to define is $H_{jjiji}^e$, which is expressed in Eq. (27).

$$H^e_{jjiji} = \frac{\partial x_{ji}}{\partial x_{j0}} \frac{\partial g^e_{ji}}{\partial x_{ji}} = 2 \frac{\partial x_{ji}}{\partial x_{j0}} \frac{\partial}{\partial x_{ji}}$$
$$\times \left( \left(1 - \frac{L_j(x_{j0})}{l_{ji}(x_{ji})}\right) \frac{\partial x_{ji}}{\partial x_{j0}} \delta_{ji} \right)$$
$$= 2 \frac{\partial x_{ji}}{\partial x_{j0}} \left( \frac{\partial x_{ji}}{\partial x_{j0}} \delta_{ji} \left( \frac{\partial}{\partial x_{ji}} \left(1 - \frac{L_j(x_{j0})}{l_{ji}(x_{ji})}\right) \right) \right)^T$$
$$+ \left(1 - \frac{L_j(x_{j0})}{l_{ji}(x_{ji})}\right) \frac{\partial}{\partial x_{ji}} \left( \frac{\partial x_{ji}}{\partial x_{j0}} \delta_{ji} \right) \right)$$
(27)

The last term of Eq. (27) is equal to zero, due to the fact that:

$$\frac{\partial x_{ji}}{\partial x_{j0}} \delta_{ji} = \begin{Bmatrix} f_k(x_{jk0} - x_{jl0}) \\ f_k(y_{jk0} - y_{jl0}) \\ -f_l(x_{jk0} - x_{jl0}) \\ -f_l(y_{jk0} - y_{jl0}) \end{Bmatrix}$$
(28)

Because if $f_k = 1$, $x_{jki}$ is not an independent variable, because it equals $x_{jk0}$. Similar reasoning can be made for the rest of the elements. So one has the expression in Eq. (29):

$$\begin{aligned} H^e_{jjiji} &= 2 \frac{\partial x_{ji}}{\partial x_{j0}} \left( \frac{\partial x_{ji}}{\partial x_{j0}} \delta_{ji} \left( \frac{\partial}{\partial x_{ji}} \left(1 - \frac{L_j(x_{j0})}{l_{ji}(x_{ji})}\right) \right)^T \right) \\ &= 2 \frac{L_j(x_{j0})}{l_j^3(x_{ji})} \frac{\partial x_{ji}}{\partial x_{j0}} \left( \delta_{ji} \delta_{ji}^T \left( \frac{\partial x_{ji}}{\partial x_{j0}} \right)^T \right) \end{aligned}$$
(29)

## 7 Optimization algorithm

As exposed before, the chosen optimization algorithm is an in-house developed SQP method, which, in our case, has full Hessian analysis. This is necessary because, as exposed before, when using initial coordinates as parameters of the synthesis, the Hessian matrix should always be underdetermined. This algorithm is applied both to the synthesis problem and the inner deformed position function, which, as exposed before, is itself an optimization problem. The Hessian matrices and gradient vectors are built via assembly of the trusses matrices and, afterwards, linear restrictions (required for the inner function) are introduced via direct manipulation of these matrices. The resultant linear system is afterwards diagonalized by means of the method presented in [30], which is able to solve underdetermined systems.

This allows one to obtain the increment vector in a decoupled system, where one can verify the signs of each variable to check if it leads to a maximization or a minimization, or an inflexion point. The underdetermined nature of the problem will lead to, at least, an stationary point in one direction. After this procedure, unidimensional optimization techniques are applied.

The optimization algorithm chosen is of the exploitative type. This means that it is very effective when one wants to improve an initial mechanism with an acceptable quality. If it is desired to find an appropriate mechanism from a fresh start, it would be more logic to use one explorative algorithm such as the Genetic Algorithm. In such case, it would not be possible to use this deformation energy based error function because, as it is demonstrated in reference [27], this type of functions leads to mechanisms with low stiffness and, therefore, of low usefulness.

## 8 Experimental results

In order to verify the behaviour of the algorithm, some simple examples are shown. The first topology that will be addressed consists on a simply articulated truss, which is wanted to follow a prescribed path. Taking into account the fact that any point of a simply articulated truss describes a circle, the result should be an adjustment to a circle as seen in Fig. 6. Here a truss is drawn with node A fixed while B describes the prescribed path defined by precision points 0–1–2.

In the left-handed picture the problem along with the starting guess is shown. The obtained result is shown in the right-handed picture, where it can be observed the optimized position of the fixed node A.

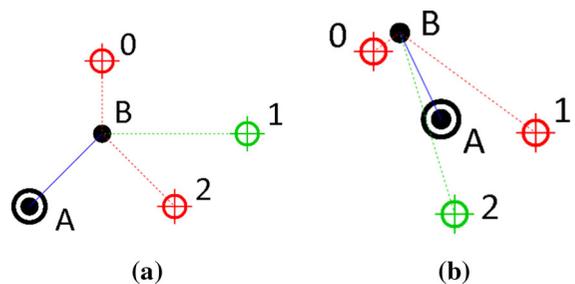

**Fig. 6** Simple truss, following a circle described by three points

The result is obtained in about 8 iterations to a precision in the order of $10^{-31}$.

Being this a numerical algorithm its sensitivity is determined by the size of the floating point used and other factors such as the preciseness of the criteria of convergence. In this case floating point of double precision have been worked with and criteria of convergence have been adjusted to obtain results with at least 5 significant numbers.

It is important to point out that although in the results presented the precision points are achieved in the specified order, this is due to the fact that they belong to feasible paths for the mechanisms of the considered typology. That is, in the optimization process it has not been introduced any condition to verify this order. However, in the proposed algorithm constraints could be introduced to force the mechanism to follow a certain order. In any case, these constraints could cause the lack of convergence towards a quality solution.

The second example deals with the same topology, but now the result is not exact. In this case it is desired to adjust the circle to 5 points as it can be appreciated in Fig. 7. While A is a fixed node, B approaches the prescribed path defined by 0–1–2–3–4.

The result is obtained in a similar number of iterations, with an increased cost for each of them. These results show that the algorithm is able to deal with both exact and approximate synthesis. In both of these examples, due to their particular nature, coupled and uncoupled formulations coincide. Now more complex problems will be addressed.

The next one is a fourbar mechanism, which is wanted to describe a 9 point path (see Fig. 8). The initial mechanism is defined by the set of coordinates expressed in Table 1. In this example fixed nodes are A

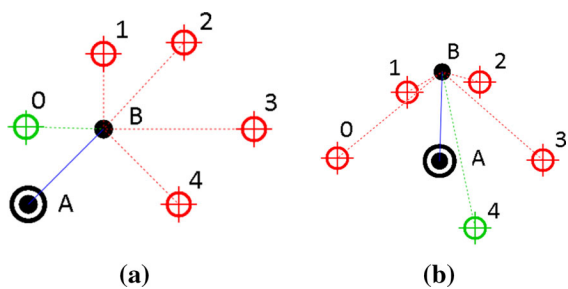

**Fig. 7** Simple truss, following a circle described by five points

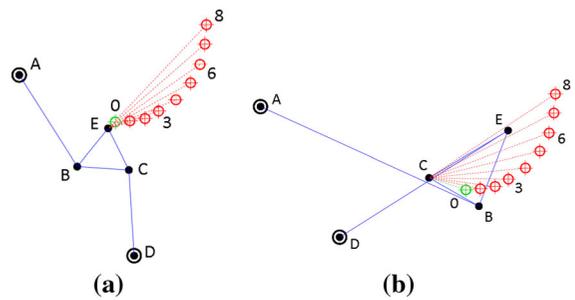

**Fig. 8** Initial position of the fourbar (left) and its obtained solution (right) with the dimensional parameters

and D, as can be seen addressed in the aforementioned figure:

And the target points are defined in Table 2.

The initial fitness is 17.2888. The algorithm reaches 0.002953. The final result is stopped due to the fact that the gradient reaches a change in the configuration, which in turn leads to an increment in the fitness for the iteration, which is not allowed in the algorithm. Some of the undeformed minimum distance points are shown in Fig. 9.

The final coordinates for the mechanism in this example are given in Table 3.

Convergence rate in the first iterations is quite good, while afterwards, the fact that it is being used an approximation of the Hessian penalizes it. Anyway, a quite good improvement is done in about 8 iterations as shown in the graphic in Fig. 10.

Obviously, with the dimensional approach, one cannot include the basement locations as optimization variables without introducing complex modifications, as explained before. In order to compare methods, the same problem will be solved including restrictions so the fixed nodes are not part of the optimization.

Using initial coordinates, the finally obtained result yields a deformation energy of 0.000615813. It may surprise that the obtained result has better fitness than the case where fixed nodes are part of the optimization, but the reason behind it is that this limitation is that the algorithm has converged to a different local minimum as shown in Fig. 11.

The minimum distance positions to the target points are quite accurate. Some of them are shown in Fig. 12.

The resultant coordinates are shown in Table 4.

When employing dimensions, the obtained result is about the same, only differenced by the convergence

**Table 1** Initial coordinates of the fourbar

| $X_A$ | $Y_A$ | $X_B$ | $Y_B$ | $X_C$ | $Y_C$ |
|---|---|---|---|---|---|
| − 5.7114 | 2.5202 | − 3.8503 | − 0.4130 | − 2.1952 | − 0.5217 |

| $X_D$ | $Y_D$ | $X_E$ | $Y_E$ |
|---|---|---|---|
| − 2.0260 | − 3.2762 | − 2.8596 | 0.8072 |

**Table 2** Coordinates of the 9 precision points to be followed by the fourbar

| X | Y | X | Y |
|---|---|---|---|
| − 2.6301 | 1.0126 | − 0.2139 | 2.2690 |
| − 2.1589 | 1.0488 | 0.0882 | 2.8610 |
| − 1.6757 | 1.1213 | 0.2443 | 3.5135 |
| − 1.2408 | 1.3630 | 0.2931 | 4.1358 |
| − 0.6850 | 1.7254 | | |

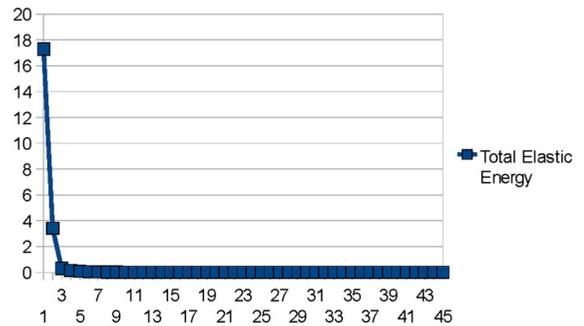

**Fig. 10** Evolution of the fitness, the elastic energy function

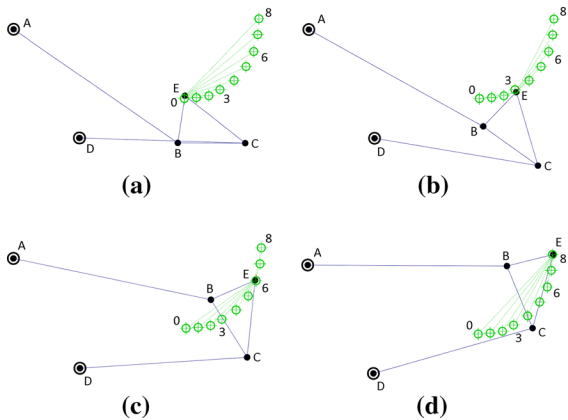

**Fig. 9** Solution mechanism in positions 0, 3, 6, and 8

**Table 3** Final coordinates of the fourbar

| $X_A$ | $Y_A$ | $X_B$ | $Y_B$ | $X_C$ | $Y_C$ |
|---|---|---|---|---|---|
| − 9.3343 | 3.7231 | − 2.2052 | 0.4771 | − 1.2526 | 2.9409 |

| $X_D$ | $Y_D$ | $X_E$ | $Y_E$ |
|---|---|---|---|
| − 6.7509 | − 0.5369 | − 3.8337 | 1.4035 |

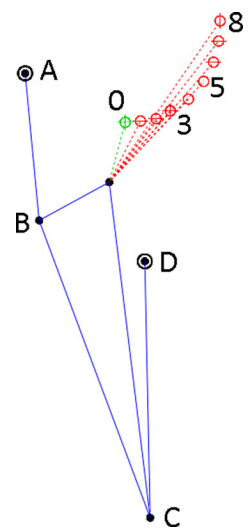

**Fig. 11** Obtained solution with the nodal coordinates parameters

moment, so one should compare the cost. The following plots in Fig. 13 describe the evolution of each of the approaches in the first 50 iterations. They are quite similar, as is the cost per iteration (although it should be slightly better in the dimensional approach, due to the reduced number of variables, but in this case the comparison is 5 to 6 and implementation differences and other factors can also affect this cost).

In this case, the resulting variables are the dimensions of the links of the mechanism, which are shown in Table 5.

As a final example, a double butterfly mechanism will be dealt with (see Fig. 14). In this case, to further show the advantages of the approach, this example will be formulated with prescribed timing: one wants the input link to achieve a determinate angle while the coupler point reaches a target position for each of the 6 precision points introduced 0–1–2–3–4–5. The A fixed node location of the input link is fixed, while the other two fixed point positions J and G are free. The initial coordinates are the ones shown in Table 6.

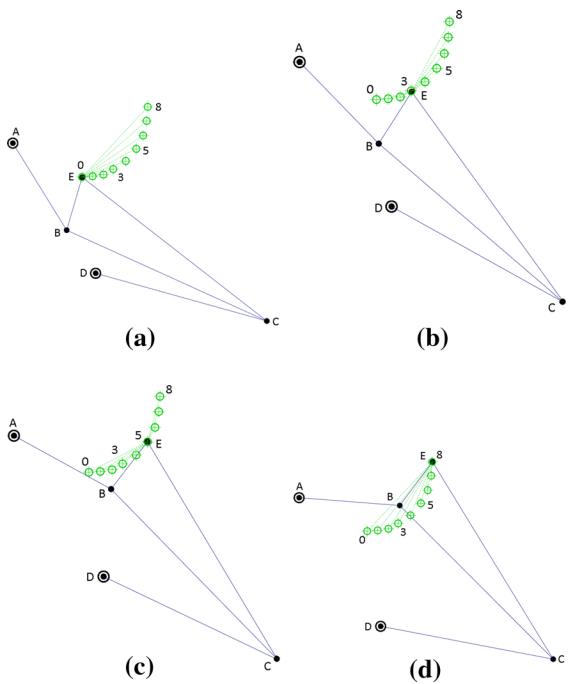

**Fig. 12** Solution mechanism in positions 0, 3, 5, and 8

The exact solution is impossible to achieve, due to the fact that there are more restrictions than variables. The obtained solution is shown in Fig. 15.

The minimum distance position to the requirements is described in Fig. 16.

The evolution of the error is as usual, quite fast at the beginning while slow at the final stages, due to the Hessian approximation as one can appreciate in

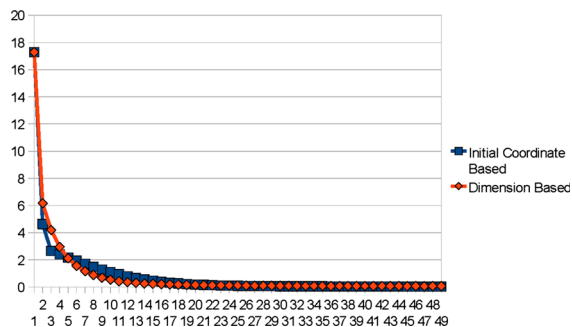

**Fig. 13** Evolution of the fitness, with synthesis based on initial coordinates (squares) and based on dimensions (diamonds)

**Table 5** Final dimensions of the fourbar

| $L_0$ | $L_1$ | $L_2$ | $L_3$ | $L_4$ |
|---|---|---|---|---|
| 2,2565 | 6,1263 | 5,4331 | 3,0782 | 4,1256 |

Fig. 17. The final solution yields 0.007, which is less than one hundred percent of the starting value.

The coordinates of the nodes in the final mechanism are given in Table 7.

Obviously, in order to successfully apply these techniques to complex mechanisms like the present one, the starting solution is of most importance, because of the presence of a large amount of local optima and also because the energy function favours low stiffness mechanisms, which can be useless, but can reach any condition. In the case of the coordinate based approach, it can also yield to degenerated 2 dof mechanisms if the initial solution is too far from the desired optima. As exposed in [1], the use of distance based functions along with genetic algorithms can give good initial solutions in these situations.

The examples shown in this work have been run on an Intel Xeon E5645@2,4GHz and the code was not programmed for multithread. The execution times are

**Table 4** final coordinates of the fourbar very variable, where the fourbar examples lie under one second, although in cases of slow convergence it has been reached, very exceptionally, the 10 min. In the example of the double butterfly presented the execution time was of 33 s.

| $X_A$ | $Y_A$ | $X_B$ | $Y_B$ | $X_C$ | $Y_C$ |
|---|---|---|---|---|---|
| − 5.7114 | 2.5202 | − 5.2815 | − 2.0171 | − 1.8628 | − 11.183 |
| $X_D$ | $Y_D$ | $X_E$ | $Y_E$ | | |
| − 2.0260 | − 3.2762 | − 3.1240 | − 0.8376 | | |

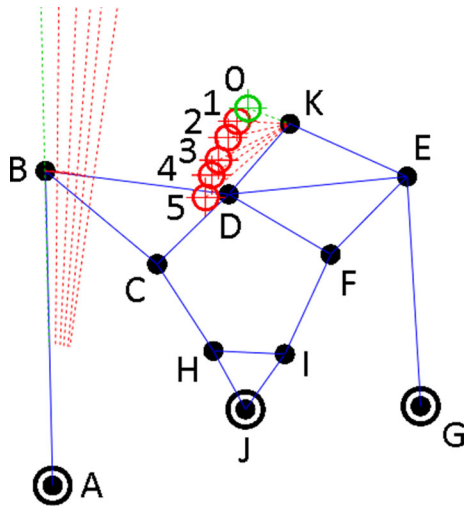

**Fig. 14** Initial guess of double butterfly with prescribed timing path generation

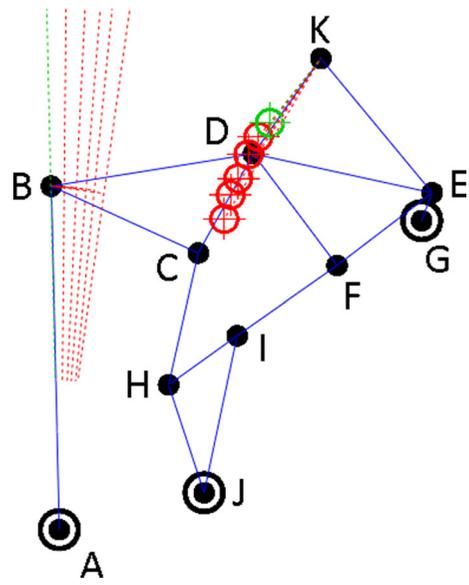

**Fig. 15** Solution of double butterfly with prescribed timing path generation

Comparing with synthesis methods based on dimensions, the use of initial coordinates presents a similar performance. This was to be expected as the number of unknowns does not increase in a considerable way.

## 9 Conclusions and future work

This paper has shown a new approach to the dimensional synthesis of mechanism which, although based in the same concepts as previous developments, introduces fundamental changes in its conception. The main contribution of this work is that thanks to the fact that the initial coordinates are used as optimization variables, the assembly configuration is included in the optimization process, which is of most impor-tance in the definition of the mechanism. A second point of interest derives from the fact that the coordinates of the fixed points are also variables of the optimization and thus, one does not need to include workarounds to optimize them. A final advantage, directly derived from the first, is that all of the possible solution vectors define a mechanism which always can be assembled, which not always holds truth when using dimensions. These advantages come to some cost, namely the fact that the same mechanism can be defined in infinite ways, thus leading to an underdetermined optimization problem. This disadvantage can successfully be overcome with an appropriate optimization method. Experimentation has shown that, depending on the problem, the use of one or another of the methods can deliver different results, so the best bet is to use both or even combinations of them. In this paper an uncoupled approach has been used, which tends to be better at the initial stages, but is slower at the final iterations. In this paper the relevant algorithms and mathematical developments have been exposed and, although they have been limited to mechanisms composed by R-Joints, they can easily be generalized to P-joints and even three-dimensional problems. In any case, the new algorithm inherits not only the advantages of the former approach, but also some of its drawbacks, specially the problem of the low stiffness mechanisms. Further developments should tackle with this problem, possibly employing a minimum distance approach, which has already shown some good results along with genetic algorithms, but requires a complex development if SQP algorithms are to be applied. The use of coupled approaches could also be of interest.

**Table 6** Initial coordinates of the double butterfly

| $X_A$ | $Y_A$ | $X_B$ | $Y_B$ | $X_C$ | $Y_C$ |
|---|---|---|---|---|---|
| − 3.7300 | − 2.0300 | − 3.8200 | 1.8900 | − 2.4300 | 0.7300 |
| $X_D$ | $Y_D$ | $X_E$ | $Y_E$ | $X_F$ | $Y_F$ |
| − 1.5400 | 1.6000 | 0.6800 | 1.8200 | − 0.2700 | 0.8500 |
| $X_G$ | $Y_G$ | $X_H$ | $Y_H$ | $X_I$ | $Y_I$ |
| 0.8500 | − 1.0400 | − 1.7300 | − 0.3500 | − 0.8400 | − 0.3900 |
| $X_J$ | $Y_J$ | $X_K$ | $Y_K$ | | |
| − 1.3300 | − 1.0800 | − 0.7800 | 2.4800 | | |

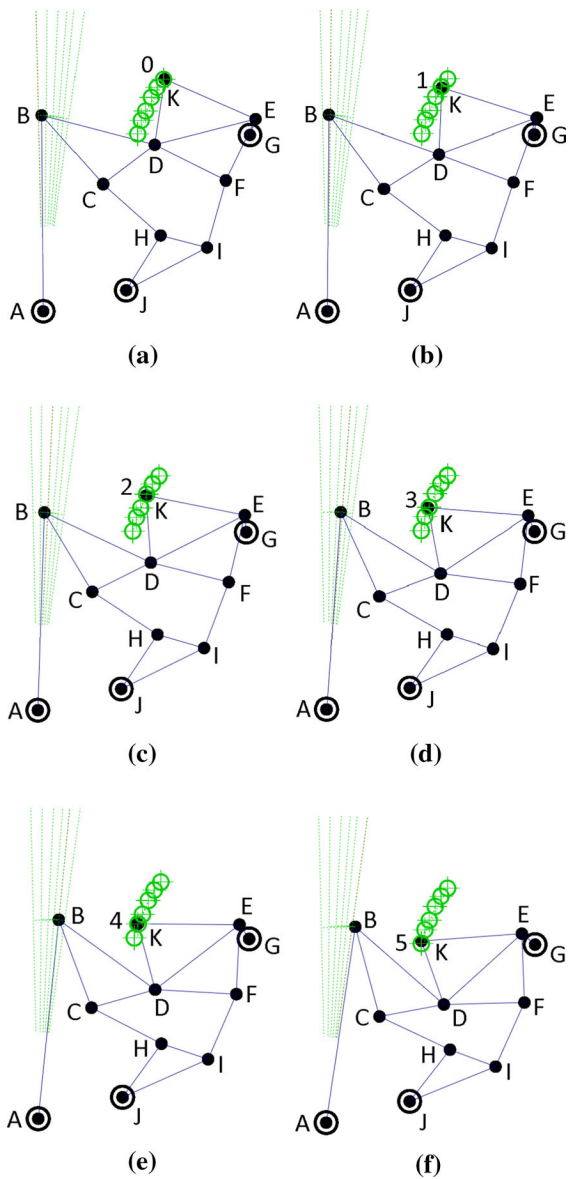

**Fig. 16** Undeformed double butterfly in minimum distance problem at every point

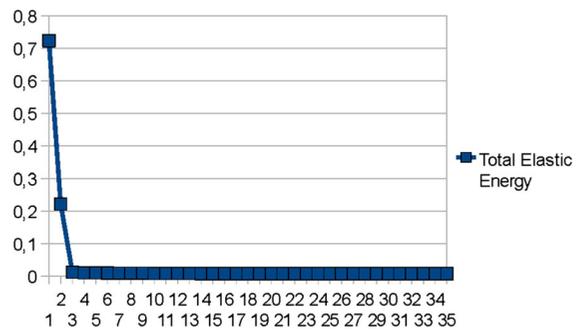

**Fig. 17** Evolution of the elastic energy function with synthesis based on initial coordinates

Table 7 Final coordinates of the double butterfly

| $X_A$ | $Y_A$ | $X_B$ | $Y_B$ | $X_C$ | $Y_C$ |
|---|---|---|---|---|---|
| − 3.7300 | − 2.0300 | − 3.8213 | 1.9397 | − 2.1265 | 1.1661 |
| $X_D$ | $Y_D$ | $X_E$ | $Y_E$ | $X_F$ | $Y_F$ |
| − 1.4903 | 2.3132 | 0.5738 | 1.8621 | − 0.5240 | 1.0219 |
| $X_G$ | $Y_G$ | $X_H$ | $Y_H$ | $X_I$ | $Y_I$ |
| 0.4528 | 1.5409 | 2.4686 | − 0.3562 | − 1.6774 | 0.2122 |
| $X_J$ | $Y_J$ | $X_K$ | $Y_K$ | | |
| − 2.0583 | − 1.6033 | − 0.7114 | 3.4163 | | |

**Acknowledgements** The authors wish to thank the Spanish Ministry of Economy and Competitiveness for its support through Grant DPI2013-46329-P and DPI2016-80372-R. Additionally the authors wish to thank the Education Department of the Basque Government for ist support through grant IT947-16.

**Compliance with ethical standards**

**Conflict of interest** The authors declare that they have no conflict of interest.